\documentclass[a4paper,11pt]{article}
    \usepackage[top=2.5cm,bottom=2.5cm,left=2.5cm,right=2.5cm]{geometry}
    \usepackage{amsmath, amssymb}
    \allowdisplaybreaks

\usepackage{mathrsfs}
\usepackage{color}
\usepackage{amsmath, epsfig, cite}
\usepackage{amssymb}
\usepackage{amsfonts}
\usepackage{latexsym}
\usepackage{graphicx}
\usepackage{enumerate}
\usepackage{enumitem}
\usepackage{booktabs}
\usepackage{overpic}
\usepackage{longtable,booktabs}
\usepackage{caption}
\usepackage[title]{appendix}

\newtheorem{df}{Definition}[section]
\newtheorem{thm}{Theorem}[section]

\newtheorem{lem}{Lemma}[section]
\newtheorem{rem}{Remark}[section]

\newcommand{\qed}{{\hfill\rule{4pt}{7pt}}}

\begin{document}
\vspace*{20mm} \pagestyle{plain}

\begin{center}
{\Large\bf Uniformly most reliable three-terminal graph of dense graphs}
\end{center}

\begin{flushleft}
\normalsize Sun Xie$^{1,2,3,4}$, Haixing Zhao$^{1,2,3}$, Jun Yin$^{1,2,3}$\\[2mm]
\footnotesize $^{1}$College of Computer, Qinghai Normal University, Xining 810008, China\\
\footnotesize $^{2}$Key Laboratory of Tibetan Information Processing, Ministry of Education, Xining 810008, China\\
\footnotesize $^{3}$Key Laboratory of Tibetan Information Processing and Machine Translation in QH, Xining 810008, China\\
\footnotesize $^{4}$Teaching and Research Department of Basic Courses, Qinghai University, Xining 810008, China\\[2mm]

Correspondence should be addressed to Haixing Zhao; h.x.zhao@163.com
\end{flushleft}

{{\bf Abstract}

\small A graph $G$ with $k$ specified target vertices in vertex set is a $k$-terminal graph. The $k$-terminal reliability is the connection probability of the fixed $k$ target vertices in a $k$-terminal graph when every edge of this graph survives independently with probability $p$. For the class of two-terminal graphs with a large number of edges, Betrand, Goff, Graves and Sun constructed a locally most reliable two-terminal graph for $p$ close to $1$, and illustrated by a counterexample that this locally most reliable graph is not the uniformly most reliable two-terminal graph. At the same time, they also determined that there is a uniformly most reliable two-terminal graph in the class obtained by deleting an edge from the complete graph with two target vertices.
This article focuses on the uniformly most reliable three-terminal graph of dense graphs with $n$ vertices and $m$ edges. First, we give the locally most reliable three-terminal graphs of $n$ and $m$ in certain ranges for $p$ close to $0$ and $1$.
Then, it is proved that there is no uniformly most reliable three-terminal graph with specific $n$ and $m$, where $n\geq7$ and $\binom{n}{2}-\lfloor\frac{n-3}{2}\rfloor\leq m\leq\binom{n}{2}-2$. Finally, some uniformly most reliable graphs are given for $n$ vertices and $m$ edges, where $4\leq n\leq 6$ and $m=\binom{n}{2}-2$ or $n\geq5$ and $m=\binom{n}{2}-1$.
}

\baselineskip 6mm

\section{Introduction}
Network reliability is a hot topic which has been general investigated using graph theoretic models. Network with $n$ vertices and $m$ edges can be modeled as a graph $G$ with the same number of vertices, edges, and interconnections as the network. For all-terminal reliability (connection probability of all vertices of a graph), lots of authors investigate the existence of a uniformly most reliable (all-terminal) graph for various values of $n$ and $m$ {\upshape\cite{Archer,Romero,Boesch,AthSobel,Brown,Gross,Myrvold,Kelmans}. However, the research on $k$-terminal reliability (connection probability of $k$ target vertices in a graph, where $2\leq k<n$) is mainly about the algorithm of computing the $k$-terminal reliability polynomial {\upshape\cite{ZhangZhao,SatyanarayanaWood,Wang,Chi,Niu}}, but only a few results on the construction of the uniformly most reliable $k$-terminal graph.

In {\upshape\cite{Kelmans}}, Kelmans has shown that for $\binom{n}{2}-\lfloor\frac{n}{2}\rfloor\leq m\leq\binom{n}{2}-2$, the uniformly most reliable graph is a complete graph with a matching removed (the matching of a graph is a set of edges with no common vertices between each other). In fact, the design of the real network often only needs to ensure the connectivity of $k$ $(2\leq k<n)$ critical vertices (target vertices). Therefore, the construction of the most reliable $k$-terminal graph has high application value. There are a few researches on the construction of the most reliable $k$-terminal structure.
In {\upshape\cite{BertrandGoffGraves}}, Betrand \emph{et al}. proved that there is no uniformly most reliable two-terminal graph when $\binom{n}{2}-\lfloor\frac{n-2}{2}\rfloor\leq m\leq\binom{n}{2}-2$. At the same time, they also proved that when $m=\binom{n}{2}-1$, the uniformly most reliable two-terminal graph is a complete graph with removing an edge between non-target vertices.
It is natural to consider the following problems.

{\bf Problems:} For the three-terminal graphs with a large number of edges, is there a uniformly most reliable graph? If it exists, what is its construction?
If it does not exist, can we construct the locally most reliable three-terminal graph?

With these questions, we further study the existence of uniformly most reliable three-terminal graphs for large $m$. It is difficult to find the exact cases that three target vertices are connected, which is NP-complete {\upshape\cite{Valiant}. In this paper, we just consider the uniformly most reliable three-terminal dense graphs.
In section $2$, some related basic definitions and notations are given.
In section $3$, the locally most reliable three-terminal graphs for given $m$ are determined. We show that there is no uniformly most reliable three-terminal graph for $n$ vertices and $m$ edges where $n\geq7$ and $\binom{n}{2}-\lfloor\frac{n-3}{2}\rfloor\leq m\leq\binom{n}{2}-2$ and give the uniformly most reliable graphs for $4\leq n\leq6$ and $m=\binom{n}{2}-2$.
In section $4$, a uniformly most reliable three-terminal graph with $n\geq5$ vertices and $m=\binom{n}{2}-1$ edges is determined.
The results are summarized in Section $5$.

\section{Basic concepts and notations}
For notations and terminologies not defined here we refer to {\upshape\cite{Bondy}}.
Let $\beta_{G}(H)$ denote the number of subgraphs which is isomorphic to $H$ in $G$. The complement of $G$, denoted by $\overline{G}$, is the graph obtained by deleting the edges of $G$ and adding edges between all non-adjacent vertices in $G$.
If $u,v\in V(G)$, then $G\cup \{uv\}$ denotes the addition of the edge $uv$ to $G$, and $G-\{uv\}$ denotes the deletion of the edge $uv$ from $G$.
The union of simple graphs $G$ and $H$, denoted by $G\cup H$, is the graph with vertex set $V(G)\cup V(H)$ and edge set $E(G)\cup E(H)$. Let $P_n$ be the path with $n$ vertices, $K_n$ be the complete graph with $n$ vertices, and $S_{n}$ be a star with $n$ vertices and $n-1$ edges.

A graph $G$ with three specified target vertices $r,s$ and $t$ in $V(G)$ is a three-terminal graph. Using $\mathcal{G}_{n,m}$ denotes the set of all simple three-terminal graphs with $n$ vertices and $m$ edges. The connectivity probability of the three specified target vertices $r,s,t$ in graph $G\in\mathcal{G}_{n,m}$ when each edge of $G$ survives independently with a fixed probability $p$ is called the three-terminal reliability of $G$ (or the three-terminal reliability polynomial of $G$), denoted by $R_3(G;p)$. A $v_1v_2\cdots v_j$-subgraph is a subgraph of $G$ in which vertices $v_1,v_2,\cdots, v_j$ are connected in the subgraph. In particular, if the $v_1v_2\cdots v_j$-subgraph with $i$ edges does not contain any $v_1v_2\cdots v_j$-subgraph with less than $i$ edges, then it is minimal, otherwise it is non-minimal. Clearly, $rst$-subgraph is a subgraph of $G$ in which three target vertices $r,s,t$ are connected in the subgraph. Similar to the definition of two-terminal reliability {\upshape\cite{BertrandGoffGraves}}, the three-terminal reliability polynomial of the graph $G\in\mathcal{G}_{n,m}$ can be written as
\begin{eqnarray}
R_3(G;p)=\sum_{i=2}^mN_i(G)p^i(1-p)^{m-i},
\end{eqnarray}
where $N_i(G)$ (or simply $N_i$) is the number of $rst$-subgraphs of graph $G$ with $i$ edges.

Similar as the definitions of the uniformly most reliable two-terminal graph {\upshape\cite{BertrandGoffGraves}} and the locally most reliable all-terminal graph {\upshape\cite{Brown}}, we defined the uniformly most reliable graph and the locally most reliable graph for three terminal graphs.

\begin{df}{\upshape}\label{df2-2}
A graph $G$ is the uniformly most reliable graph in $\mathcal{G}_{n,m}$, if $R_3(G;p)\geq R_3(H;p)$ for all $H\in\mathcal{G}_{n,m}$ and all $0\leq p\leq1$. In particular, for $p_0=0$ $($or $1 )$, if there is an $\varepsilon>0$ such that $R_3(G;p)\geq R_3(H;p)$ for all $H\in\mathcal{G}_{n,m}$ and for all $p\in[0,1]\cap(p_0-\varepsilon,p_0+\varepsilon)$, then $G$ is the locally most reliable graph in $\mathcal{G}_{n,m}$ for $p$ close to $0$ $($or for $p$ close to $1 )$.
\end{df}

\noindent
{\bf Example 1.} Figure $1$ shows all types of simple three-terminal graph in $\mathcal{G}_{4,4}$ with three target vertices $r,s,t$. Each edge of these graphs survives independently with probability $p$.

\begin{figure}[!hbpt]
\begin{center}
\includegraphics[scale=0.61]{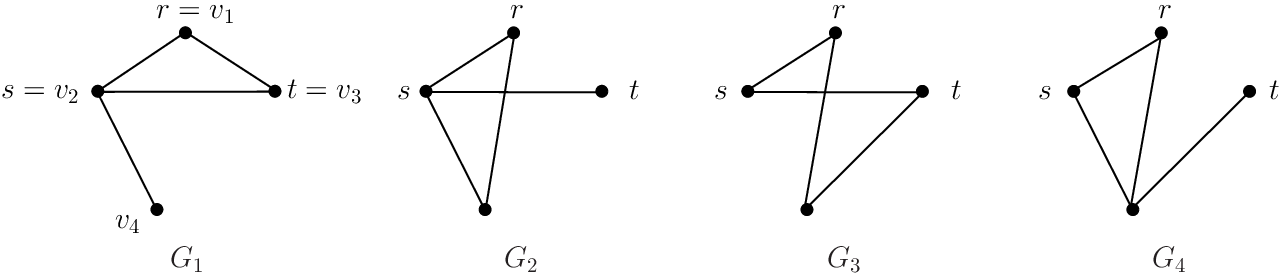}\\ [0.1cm]
Figure 1: All simple three-terminal graphs in $\mathcal{G}_{4,4}$ with three target vertices $r,s,t$.
\end{center}\label{fig20}
\end{figure}

In $G_1$, the $rst$-subgraphs with 2 edges are $\{rs,rt\}$, $\{rs,st\}$,$\{rt,st\}$,
and the $rst$-subgraphs with 3 edges are $\{rs,rt,st\}$, $\{rs,rt,sv_4\}$, $\{rs,st,sv_4\}$, $\{rt,st,sv_4\}$, while the
$rst$-subgraph with 4 edges is $\{rs,rt,st,sv_4\}$.
Obviously, $N_2(G_1)=3$, $N_3(G_1)=4$ and $N_4(G_1)=1$.
Similarly, we can calculate $N_i(G_j)$, $2\leq i, j\leq4$, which are $N_2(G_2)=1$, $N_3(G_2)=3$, $N_4(G_2)=1$; $N_2(G_3)=1$, $N_3(G_3)=4$, $N_4(G_3)=1$; $N_2(G_4)=0$, $N_3(G_4)=3$ and $N_4(G_4)=1$.
Figure $2$ shows a visualization of reliability polynomials for all graphs in $\mathcal{G}_{4,4}$. Clearly, for all $0<p<1$, $R_3(G_1;p)>R_3(G_3;p)>R_3(G_2;p)>R_3(G_4;p)$, so, $G_1$ is the uniformly most reliable graph in $\mathcal{G}_{4,4}$.

\begin{figure}[!hbpt]
\begin{center}
\includegraphics[scale=0.2]{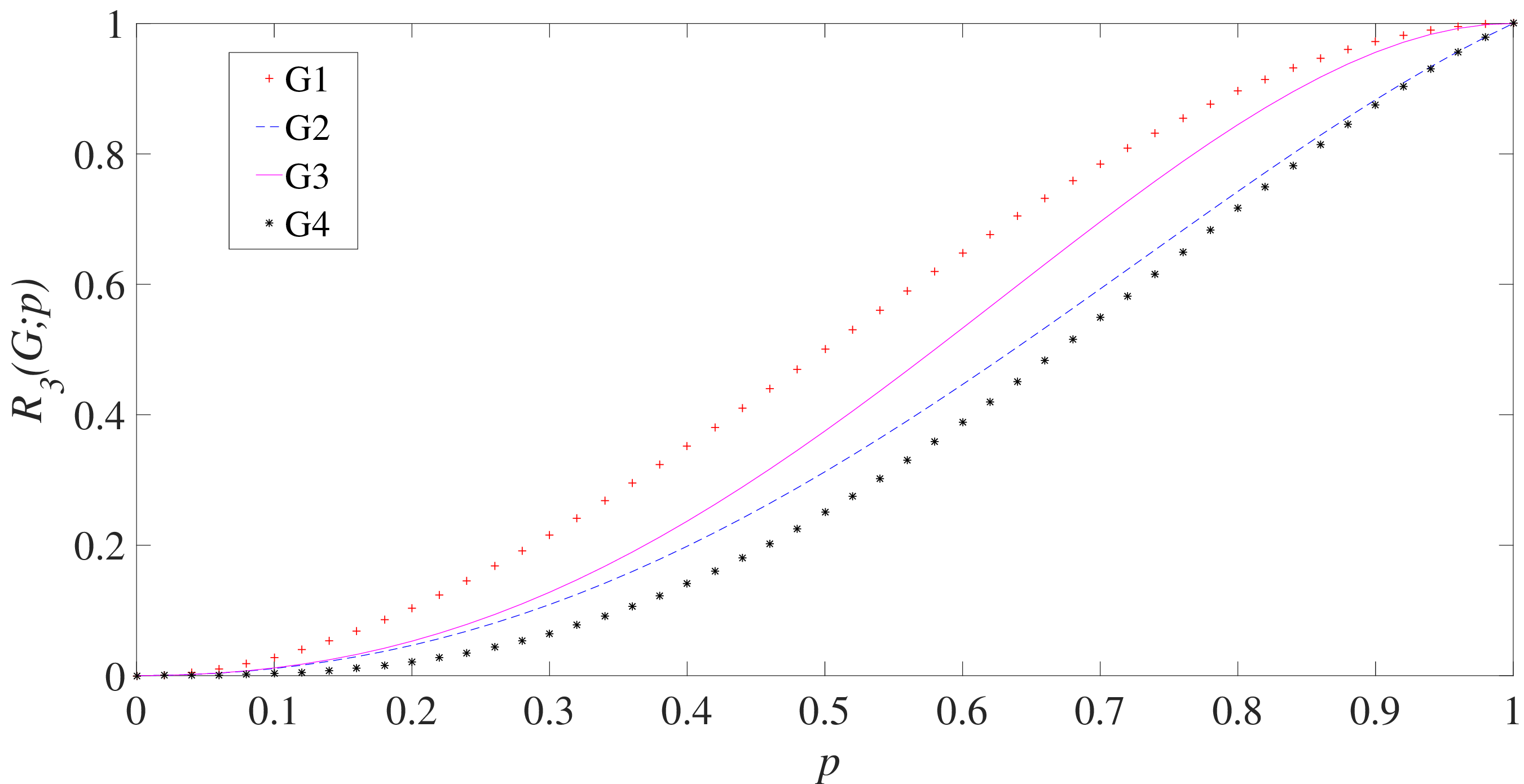} \\ [0.1cm]
Figure 2: A visualization of reliability polynomials for graphs in $\mathcal{G}_{4,4}$.
\end{center}\label{fig20}
\end{figure}

\noindent
{\bf Example 2.}  Figure $3$ shows two special simple three-terminal graphs in $\mathcal{G}_{8,26}$ with three target vertices $r,s,t$. Each edge of these graphs survives independently with probability $p$.
By calculation, we give a plot of $R_3(H_1;p)-R_3(H_2;p)$ as shown in Figure $4$. Clearly, $R_3(H_1;p)>R_3(H_2;p)$ for $p$ close to 1, and $R_3(H_1;p)<R_3(H_2;p)$ for $p$ close to 0. This article later proves that $H_1$ is the locally most reliable graph for $p$ close to 1 and $H_2$ is the locally most reliable graph for $p$ close to 0 in $\mathcal{G}_{8,26}$.

\begin{figure}[!hbpt]
\begin{center}
\includegraphics[scale=0.52]{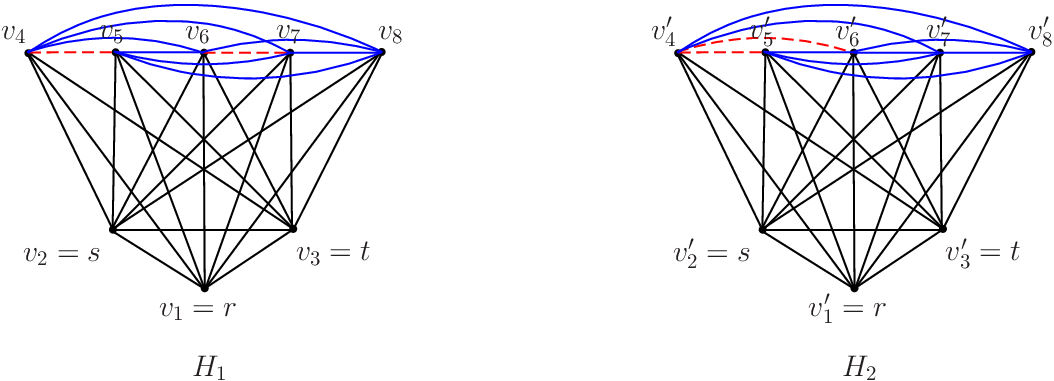} \\ [0.1cm]
Figure 3. Two special three-terminal graphs in $\mathcal{G}_{8,26}$ with three target\\
\qquad \qquad vertices $r,s,t$. The red dotted lines indicate the deleted edges.
\end{center}\label{fig20}
\end{figure}

\begin{figure}[!hbpt]
\begin{center}
\includegraphics[scale=0.18]{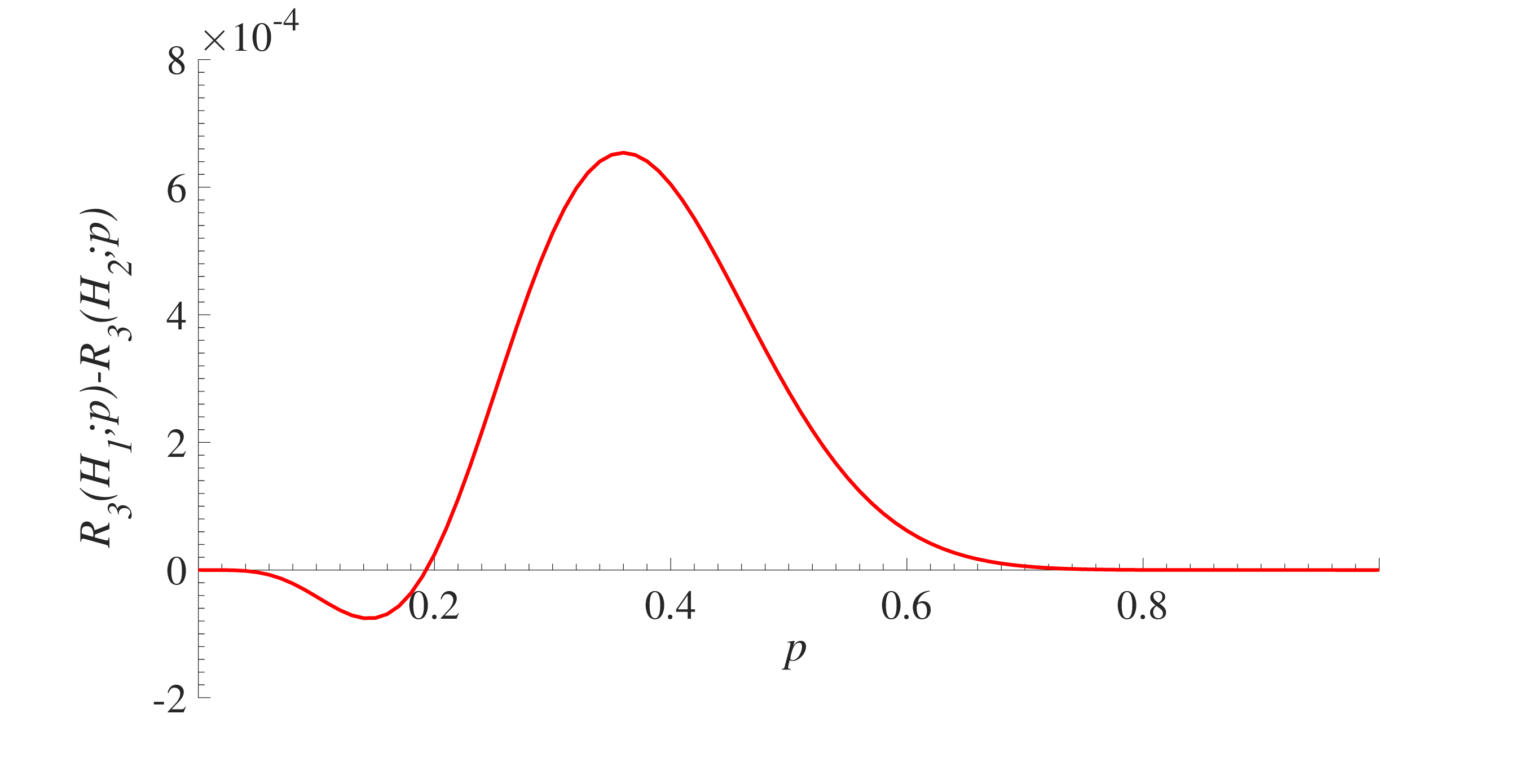} \\ [0.1cm]
Figure 4. A plot of $R_3(H_1;p)-R_3(H_2;p)$.
\end{center}\label{fig20}
\end{figure}

Many researches focus on determining a uniformly most reliable graph for given number of vertices $n$ and edges $m$, as shown in {\upshape\cite{Archer,BertrandGoffGraves,Romero}}. If there is no uniformly most reliable graph, researchers usually focus on determining the locally most reliable graph for $p$ close to $0$ or $1$, as shown in {\upshape\cite{Byer,Brown}}. At present, there are few studies on determining whether there is a uniformly most reliable three-terminal graph. Therefore, in this paper, we study the uniformly most reliable graph and the locally most reliable graph of dense three-terminal graph. One can see that some dense graphs have uniformly most reliable three-terminal graphs as Example $1$, some does not have the uniformly most reliable graph, but has locally most reliable three-terminal graphs as Example $2$.

\section{Some locally most reliable three-terminal graphs}
In this section, the locally most reliable three-terminal graph with $n\geq7$ and $m\in[\binom{n}{2}-(n-4),\binom{n}{2}-2]$ for $p$ close to 0 is determined and the locally most reliable three-terminal graph with $n\geq7$ and $m\in[\binom{n}{2}-\lfloor\frac{n-3}{2}\rfloor,\binom{n}{2}-2]$ for $p$ close to 1 is also determined. Then, it is shown that for $n\geq7$ and $m\in[\binom{n}{2}-\lfloor\frac{n-3}{2}\rfloor,\binom{n}{2}-2]$, there is no uniformly most reliable graph in $\mathcal{G}_{n,m}$. For $n=4,5,6$ and $m=\binom{n}{2}-2$, there is a uniformly most reliable three-terminal graph. To prove these results, we first introduce some related lemmas.

\begin{lem}{\upshape \cite{Byer}}\label{lem3-1}
Let $n\geq1$ and $0\leq m\leq n-1$ be positive integers.
If $m\neq3$, then the unique simple graph on $n$ vertices and $m$ edges with the maximum number of $P_3$ is $S_{m+1}\cup \overline{K_{n-m-1}}$.
If $m=3$, then there are two simple graphs with the maximum number of $P_3$: $K_3\cup \overline{K_{n-3}}$ and $S_{4}\cup \overline{K_{n-4}}$.
\end{lem}

In general, it is difficult to calculate the three-terminal reliability polynomial of graph. Therefore, we study the locally most reliable graph by the following Lemma 3.2, which is extracted from {\upshape\cite{BertrandGoffGraves}.

\begin{lem}\label{lem3-2}
Let the three-terminal reliable polynomials of $G,H\in\mathcal{G}_{n,m}$ be

$R_3(G;p)=\sum\limits_{i=2}^mN_i(G)p^i(1-p)^{m-i}$ \quad and \quad $R_3(H;p)=\sum\limits_{i=2}^mN_i(H)p^i(1-p)^{m-i}.$

\noindent
Let $N_i(G)=N_i(H)$ for $1\leq i< k$ and for $l< i\leq m$, where $k$ and $l$ are integers. Then
\begin{enumerate}[itemindent=1em]
\item[$(1)$] For $p$ close to $0$, if $N_k(G)>N_k(H)$, then $R_3(G;p)>R_3(H;p)$,
\item[$(2)$] For $p$ close to $1$, if $N_l(G)>N_l(H)$, then $R_3(G;p)>R_3(H;p)$.
\end{enumerate}
\end{lem}

By Lemma $3.2$,  we get the following conclusions.

$(1)$ If $G\in\mathcal{G}_{n,m}$ is the locally most reliable graph for $p$ close to $0$, then it must contain the triangle $rst$ and the value of $N_3(G)$ is the maximum among graphs containing the triangle $rst$ in $\mathcal{G}_{n,m}$.

$(2)$ An $rst$-cutset of $G$ is a set of edges whose deletion results in the disconnection of vertices $r,s$ and $t$ in $G$ and the number of edges is its size. The edge connectivity of $r,s$ and $t$ in $G$ is the smallest size of an $rst$-cutset of $G$, denoted by $\lambda_{rst}(G)$ or simply $\lambda(G)$. If $G\in\mathcal{G}_{n,m}$ is the locally most reliable graph for $p$ close to $1$, then it must have the maximum $\lambda_{rst}(G)$. Since $N_i=\binom{m}{i}$ ($m-\lambda+1\leq i\leq m$) and $N_{m-\lambda}=\binom{m}{\lambda}-a$, where $a$ is the number of the $rst$-cutsets with size $\lambda$, $G$ must have the minimum $a$ among the graphs with the largest edge connectivity of $r,s$ and $t$ in $\mathcal{G}_{n,m}$.

Now, we demonstrate the locally most reliable graph for three-terminal graphs for $p$ close to 0 or 1. We first introduce two important graphs for this section, as follows.

Let $n\geq7$ and $2\leq l\leq n-4$ be positive integers.
Using $A_{n,l}$ denotes the three-terminal graph on $n$ vertices and $\binom{n}{2}-l$ edges with vertex set
$V(A_{n,l})=\{r=v_1,s=v_2,t=v_3,v_4,\cdots,v_n\}$
and edge set
$E(A_{n,l})=\{v_iv_j|1\leq i<j\leq n\}-\{v_{4}v_{i+3}|2\leq i\leq l+1\}.$

Let $n\geq7$ and $2\leq l\leq\lfloor\frac{n-3}{2}\rfloor$ be positive integers.
Using $A'_{n,l}$ denotes the three-terminal graph on $n$ vertices and $\binom{n}{2}-l$ edges with vertex set
$V(A'_{n,l})=\{r=v_1,s=v_2,t=v_3,v_4,\cdots,v_n\}$
and edge set
$E(A'_{n,l})=\{v_iv_j|1\leq i<j\leq n\}-\{v_{2i}v_{2i+1}|2\leq i\leq l+1\}.$

Figure $5$ depicts these two three-terminal graphs with $11$ vertices and $51$ edges.

\begin{figure}[!hbpt]
\begin{center}
\includegraphics[scale=0.56]{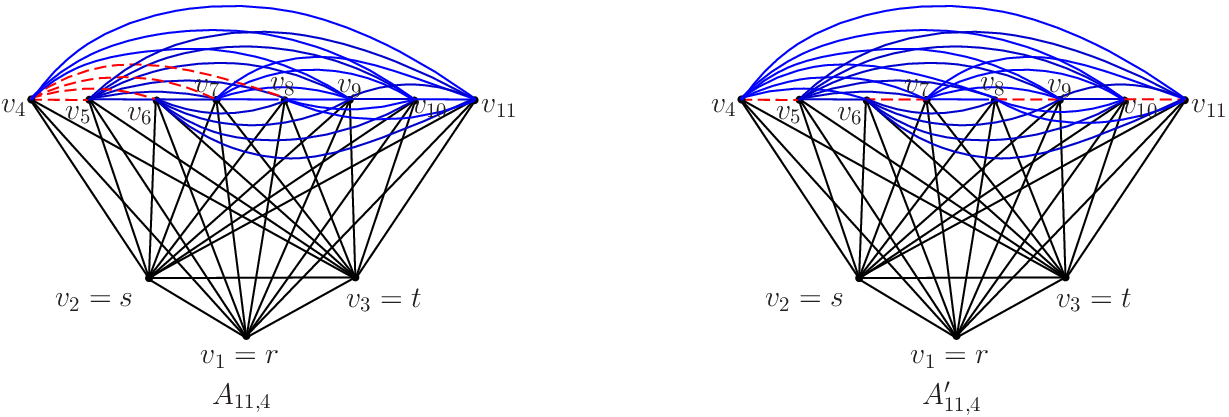} \\ [0.05cm]
Figure 5: Graph $A_{11,4}$ (left) and Graph $A'_{11,4}$ (right)\\
The red dotted lines indicate the deleted edges.
\end{center}\label{fig11}
\end{figure}

\begin{thm}\label{th3-1}
Let $n\geq7$, $2\leq l\leq n-4$ and $m=\binom{n}{2}-l$ be positive integers. Then
\begin{enumerate}
\item[$(1)$] If $l=3$, then the graph $A^*_{n,3}=A_{n,3}\cup\{v_4v_7\}-\{v_5v_6\}$ is the unique locally most reliable graph in $\mathcal{G}_{n,m}$ for $p$ close to $0$,
\item[$(2)$] If $l\neq3$, then the graph $A_{n,l}$ is the unique locally most reliable graph in $\mathcal{G}_{n,m}$ for $p$ close to $0$.
\end{enumerate}
\end{thm}

\begin{pf}
Suppose that $G$ is the locally most reliable graph in $\mathcal{G}_{n,m}$ for $p$ close to $0$. Then by Lemma $3.2$, $G$ must contain the triangle $rst$ and $N_3$ (the number of $rst$-subgraphs with 3 edges in $G$) must be maximum among graphs containing the triangle $rst$ in $\mathcal{G}_{n,m}$.

For the sake of calculating $N_3$, using $a,b,c$ denotes the number of $rst$-subgraphs with 3 edges containing 2,1,0 edges in triangle $rst$ respectively.
Then $a$ is the number of sets $\{rs,st,v_iv_j\}$, $\{rs,rt,v_iv_j\}$ and $\{rt,st,v_iv_j\}$ ($1\leq i,j\leq n$), $b$ is the number of sets $\{rs,rv_i,tv_i\}$, $\{rs,sv_i,tv_i\}$, $\{rt,rv_i,sv_i\}$, $\{rt,tv_i,sv_i\}$, $\{st,rv_i, sv_i\}$, $\{st,rv_i,tv_i\}$ ($4\leq i\leq n$) and $c$ is the number of set $\{v_ir, v_is, v_it\}$ ($4\leq i\leq n$). It is no hard to see that $N_3=a+b+c$.
Clearly, for graphs containing the triangle $rst$ in $\mathcal{G}_{n,m}$, $a$ is a constant which equals $3(m-3)+1$, and $N_3$ takes the maximum value if and only if $b$ and $c$ attains the maximum value.
Note that if $c$ takes the maximum value $n-3$, then the value of $b$ also reaches its maximum, that is, $E(G)$ contains the edges $v_ir, v_is, v_it$ for all $4\leq i\leq n$.
Now, consider the remaining $\binom{n}{2}-l-3n+6$ edges among non-target vertices in $G$ that has not been described. Since $G$ is a dense graph, it is often easier to consider the deleted $l$ edges among non-target vertices.

By Lemma $3.2$, we continue to calculate the coefficients $N_i=b_i+c_i$ ($i>3$), where $b_i$ and $c_i$ is the number of minimal $rst$-subgraphs and non-minimal $rst$-subgraphs with $i$ edges, respectively. We now compute $N_4$.
Clearly, $b_4$ is the sum of the numbers of sets $\{sv_i,v_it,sv_j,v_jr\}$, $\{sv_i,v_iv_j,v_jt,sr\}$ and $\{sv_i,v_iv_j,v_jt,v_jr\}$ ($4\leq i,j\leq n$). The non-minimal $rst$-subgraph with $4$ edges includes two cases: one is containing the minimal $rst$-subgraph with 2 edges and
the other is containing the minimal $rst$-subgraph with 3 edges but no minimal $rst$-subgraph with 2 edges.
By calculation, $b_4=6\binom{n-3}{2}+12(m-3n+6)+6(m-3n+6)$ and $c_4=3\binom{m-3}{2}+(m-3)+(n-3)(m-3)+6(n-3)(m-6)$ which are constants.
So, $N_4$ is a constant.

Furthermore, we need to calculate $N_5$.
Clearly, $b_5$ is the sum of the numbers of sets $\{sv_i,v_iv_j,v_jr,rv_k,v_kt\}$, $\{sv_i,v_iv_j$, $v_jv_k,v_kt,rt\}$, $\{sv_i,v_iv_j,v_jr,v_jv_k,v_kt\}$ and $\{rv_i,sv_i,v_iv_j,v_jv_k,v_kt\}$ ($4\leq i,j,k\leq n$). The non-minimal $rst$-subgraph with $5$ edges includes three cases:
containing the minimal $rst$-subgraph with 2 edges,
containing the minimal $rst$-subgraph with 3 edges but no minimal $rst$-subgraph with 2 edges,
containing the minimal $rst$-subgraph with 4 edges but no minimal $rst$-subgraph with less than 4 edges.
By calculation, $c_5=3\binom{m-3}{3}+\binom{m-3}{2}+7(n-3)\binom{m-6}{2}+3(n-3)(m-6)-12\binom{n-3}{2}+18(m-3n+6)(m-10)+6(m-3n+6)+6(m-9)\binom{n-3}{2}$, which is a constant. And
\begin{equation}
\begin{split}
b_5&= 12(m-3n+6)(n-5)+24\sum_{i=4}^n\binom{d(v_i)-3}{2}\\
&= 12\sum_{i=4}^n(d(v_i)-3)^2+12(m-3n+6)(n-7)\\
&= 24\sum_{i=4}^n\binom{d_i}{2}+24l+12(n-3)(n-4)^2-48(n-4)l+12(m-3n+6)(n-7),
\end{split}
\end{equation}

where $d(v_i)$ is the degree of $v_i$ in $G$, and $d_i$ is the number of edges deleted on the non-target vertex $v_i$ $(4\leq i\leq n)$, which is the degree of $v_i$ in $\overline{G}$.
Note that the number of subgraphs as $P_3$ in $\overline{G}$ is $\sum\limits_{i=4}^n\binom{d_i}{2}$. Obviously, $N_5$ of $R_3(G;p)$ is related to $\beta_{\overline{G}}(P_{3})$, the number of subgraphs as $P_3$ in $\overline{G}$. There are at least three isolate vertices in $\overline{G}$, which are target vertices in $G$. $\beta_{\overline{G}}(P_{3})$ is corresponding with $n-3$ non-target vertices.

$(1)$ By Lemma $3.1$, if $l=3$, then the number of $P_3$ in a simple graph with $n-3$ vertices and $l$ edges reaches the maximum if the graph is either $K_3\cup \overline{K_{n-6}}$ or $S_{4} \cup \overline{K_{n-7}}$. Thus, $\overline{G}$ is either $K_3\cup \overline{K_{n-6}}\cup\{r,s,t\}$ or $S_{4} \cup \overline{K_{n-7}}\cup\{r,s,t\}$, which means that $G$ is either $A_{n,3}$ or $A^*_{n,3}$, where $N_5(A_{n,3})=N_5(A^{*}_{n,3})$.
So, we need to compare $N_6(A_{n,3})$ and $N_6(A^{*}_{n,3})$, which can be calculated as the same analysis of $N_5$. For convenient, if $H\in\mathcal{G}(n,m)$, then let $\widehat{H}$ be the graph $H-\{s,t\}$. Then
$N_6(A_{n,3})=\binom{m-3}{3}+3\binom{m-3}{4}+3(n-3)\binom{m-6}{2}+7(n-3)\binom{m-6}{3}
+[6\binom{m-9}{2}-12(m-9)-13]\binom{n-3}{2}-30\binom{n-3}{3}+[18\binom{m-10}{2}+6 (m-10)-36(n-5)+12(n-5)(m-13)+12(m-3n+7)](m-3n+6)
+(12n+24m-459)\beta_{\widehat{A_{n,3}}}(P_3)-66\beta_{\widehat{A_{n,3}}}(K_3)+30\beta_{\widehat{A_{n,3}}}(P_4)+6\beta_{\widehat{A_{n,3}}}(S_{4})$ and
$N_6(A^{*}_{n,3})=\binom{m-3}{3}+3\binom{m-3}{4}+3(n-3)\binom{m-6}{2}+7(n-3)\binom{m-6}{3}
+[6\binom{m-9}{2}-12(m-9)-13]\binom{n-3}{2}-30\binom{n-3}{3}+[18\binom{m-10}{2}+6 (m-10)-36(n-5)+12(n-5)(m-13)+12(m-3n+7)](m-3n+6)+
(12n+24m-459)\beta_{\widehat{A^{*}_{n,3}}}(P_3)-66\beta_{\widehat{A^{*}_{n,3}}}(K_3)+30\beta_{\widehat{A^{*}_{n,3}}}(P_4)+6\beta_{\widehat{A^{*}_{n,3}}}(S_{4})$,
By calculation, $\beta_{\widehat{A_{n,3}}}(P_4)=\beta_{\widehat{A^*_{n,3}}}(P_4)$, $\beta_{\widehat{A_{n,3}}}(S_{4})-\beta_{\widehat{A^*_{n,3}}}(S_{4})=-1$ and $\beta_{\widehat{A_{n,3}}}(K_3)-\beta_{\widehat{A^*_{n,3}}}(K_3)=1$.
So,
\begin{equation}
\begin{split}
N_6(A_{n,3})-N_6(A^*_{n,3})&= 30[\beta_{\widehat{A_{n,3}}}(P_4)-\beta_{\widehat{A^*_{n,3}}}(P_4)]
+6[\beta_{\widehat{A_{n,3}}}(S_{4})-\beta_{\widehat{A^*_{n,3}}}(S_{4})]\\
&\quad -66[\beta_{\widehat{A_{n,3}}}(K_3)-\beta_{\widehat{A^*_{n,3}}}(K_3)]\\
&=-72<0.
\end{split}
\end{equation}

Therefore, if $l=3$, the graph $A^*_{n,3}$ is the unique locally most reliable graph in $\mathcal{G}_{n,m}$ for $p$ close to $0$.

$(2)$ By Lemma $3.1$, if $l\neq3$, then the number of $P_3$ in a simple graph with $n-3$ vertices and $l$ edges is maximized only if the graph is $S_{l+1}\cup \overline{K_{n-4-l}}$. Thus, $\overline{G}$ is $S_{l+1}\cup \overline{K_{n-4-l}}\cup\{r,s,t\}$, which means that $G$ is $A_{n,l}$.
Therefore, if $l\neq3$, then the graph $A_{n,l}$ is the unique locally most reliable graph in $\mathcal{G}_{n,m}$ for $p$ close to $0$.\qed
\end{pf}

\begin{thm}\label{th3-2}
Let $n\geq7$, $2\leq l\leq\lfloor\frac{n-3}{2}\rfloor$ and $m=\binom{n}{2}-l$ be positive integers. Then $A'_{n,l}$ is the unique locally most reliable graph in $\mathcal{G}_{n,m}$ for $p$ close to $1$.
\end{thm}

\begin{pf}
Let $G\in\mathcal{G}_{n,m}$ be the most reliable graph for $p$ close to 1. By Lemma $3.2$, $G$ must have the largest edge connectivity of $r,s$ and $t$, which is the size of the smallest $rst$-cutset of $G$ as large as possible. For convenient, let $\mathcal{H}_{n,m}$ be the set of three-terminal graph with the largest edge connectivity of $r,s$ and $t$. Obviously, $G\in\mathcal{H}_{n,m}$.

For any graph $H\in\mathcal{H}_{n,m}$, let $C$ be the minimal $rst$-cutset of $H$, then there must exist a component containing just one target vertex and $k$ $(0\leq k\leq n-3)$ non-target vertices $u_i$ $(1\leq i\leq k)$ in $H-C$, without loss of generality, setting this target vertex as $r$.
Clearly, $\lambda(H)\leq min\{d(r),d(s),d(t)\}\leq n-1$. If $d(r)=d(s)=d(t)=n-1$, then $|C|\geq d(r)-k+2k=d(r)+k\geq n-1$. Hence, $\lambda$ can arrive at the maximum value $n-1$ if and only if $d(r)=d(s)=d(t)=n-1$. Thus, three target vertices in $H$ are adjacent to every vertex in $V(H)$. However, the connections among non-target vertices are various. In order to get the precise construction of $G$, by Lemma $3.2$, we need to consider the coefficient $N_{m-i}$ $(i\geq\lambda)$. Obviously, $N_{m-i}=\binom{m}{m-i}-m_i$, where $m_i$ is the number of the $rst$-cutsets of size $i$. For given $i$, $\binom{m}{m-i}$ is a constant and $N_{m-i}$ is corresponding with $m_i$.
Let $\delta$ be the minimum degree of $H$. Since $n\geq7$ and $\binom{n}{2}-\lfloor\frac{n-3}{2}\rfloor\leq m\leq\binom{n}{2}-2$, we have $\delta\geq n-1-\lfloor\frac{n-3}{2}\rfloor\geq4$.

Continue to calculate the minimal $rst$-cutset of $H$, $|C|=\lambda+\sum\limits_{i=1}^kd(u_i)-2k-2m'\geq \lambda+\delta k-2k-2\binom{k}{2}=\lambda-k^2+(\delta-1)k$, where $m'$ is the edge number among non-target vertices $u_1,u_2,\cdots,u_k$. If $k>\delta-2$, then $|C|\geq d(r)+k>\lambda+\delta-2$. If $1\leq k\leq\delta-2$, then $|C|\geq\lambda+\delta-2$, where the equation holds if the component containing $r$ in $H-C$ are following two cases: one is composed of target vertex $r$ and one non-target vertex with degree $\delta$, the other is composed of target vertex $r$ and $\delta-2$ non-target vertices with degree $\delta$ and the induced subgraph by vertices in this component is $K_{\delta-1}$.
So, the number of $C$ with size $\lambda+1,\lambda+2,...,\lambda+\delta-3$ is $0$. Meanwhile, the number of $C$ with size $\lambda$ is $3$. By calculation, we see that $m_{j}(H)=3\binom{m-\lambda}{j-\lambda}$ for $\lambda\leq j\leq\lambda+\delta-3$
and $m_{\lambda+\delta-2}(H)=3\binom{m-\lambda}{\delta-2}+3n_1+3n_2$, where $n_1$ is the number of non-target vertices with degree $\delta$ and $n_2$ is the number of $K_{\delta-2}$ in the induced subgraph of non-target vertices with degree $\delta$.

Similarly as the calculation of $H$, we have $m_{j}(G)=3\binom{m-\lambda}{j-\lambda}$ for $\lambda\leq j\leq\lambda+\delta(G)-3$
and $m_{\lambda+\delta(G)-2}(G)=3\binom{m-\lambda}{\delta(G)-2}+3n_1'+3n_2'$, where $n_1'$ is the number of non-target vertices with degree $\delta(G)$ and $n_2'$ is the number of $K_{\delta(G)-2}$ in the induced subgraph of non-target vertices with degree $\delta$ in $G$.
Obviously, $n_1'>0$, $n_2'\geq0$ and if $\delta(G)<\delta$, then $m_{j}(G)=m_{j}(H)$ for $\lambda\leq j\leq\lambda+\delta(G)-3$, and $m_{\lambda+\delta(G)-2}(G)>3\binom{m-\lambda}{\delta(G)-2}$ but $m_{\lambda+\delta(G)-2}(H)=3\binom{m-\lambda}{\delta(G)-2}$. So,
$N_{m-j}(G)=N_{m-j}(H)$ for $\lambda\leq j\leq\lambda+\delta(G)-3$ and $N_{m-(\lambda+\delta(G)-2)}(G)<N_{m-(\lambda+\delta(G)-2)}(H)$, which contradicts the assumption that $G$ is the most reliable graph for $p$ close to 1. Then we have $\delta(G)\geq\delta$.

Clearly, $\delta(G)=\lfloor\frac{2(m-3n+6)}{n-3}\rfloor+3$. Since $\binom{n}{2}-\lfloor\frac{n-3}{2}\rfloor\leq m\leq\binom{n}{2}-2$, the minimum degree of graph $G$ is larger than $n-2$, which implies that for each $v\in V(G)$, $d(v)\geq n-2$.

Therefore, $A'_{n,l}$ is the unique locally most reliable graph in $\mathcal{G}_{n,m}$ for $p$ close to $1$.\qed
\end{pf}

As a straightforward consequence of Theorems $3.1$ and $3.2$, we obtain the following theorem.

\begin{thm}\label{th3-3}
Let $n\geq7$ and $2\leq l\leq\lfloor\frac{n-3}{2}\rfloor$ be positive integers. If $m=\binom{n}{2}-l$, then there is no uniformly most reliable three-terminal graph in $\mathcal{G}_{n,m}$.
\end{thm}

Three specific classes of graphs with order $n$ and size $m=\binom{n}{2}-2$ are uniformly most reliable graphs.
\begin{rem}
For $n=4$ and $m=4$, there is a uniformly most reliable three-terminal graph in $\mathcal{G}_{4,4}$ $($see Example 1$)$. For $n=5$ and $m=8$, there is a uniformly most reliable three-terminal graph in $\mathcal{G}_{5,8}$ $($ graphs are shown in Figure 6 in Appendix A, and the comparison of the reliability polynomial is shown in Table $1$ in Appendix A $)$. For $n=6$ and $m=13$, there is a uniformly most reliable three-terminal graph in $\mathcal{G}_{6,13}$ $($graphs are shown in Figure 7 in Appendix A, and the comparison of the reliability polynomial is shown in Table $2$ in Appendix A $)$.
\end{rem}

\section{The uniformly most reliable three-terminal graph}
For the three-terminal graph with $\binom{n}{2}$ edges, it is complete graph, which is also the uniformly most reliable graph. In this section, we determine the uniformly most reliable graph in $\mathcal{G}_{n,m}$ with $m=\binom{n}{2}-1$.

When we remove one edge from $K_n$ with three target vertices, there are only three cases: the edge between target vertices; the edge between a target vertex and a non-target vertex; the edge between non-target vertices. Let $n\geq5$ and $m=\binom{n}{2}-1$ be positive integers.

$(1)$ Using $X_n$ denotes the three-terminal graph on $n$ vertices and $m$ edges with vertex set $V(X_n)=\{r=x_1,s=x_2,t=x_3,x_4,\cdots,x_n\}$ and edge set $E(X_n)=\{x_ix_j|1\leq i<j\leq n\}-\{rs\}$.

$(2)$ Using $Y_n$ denotes the three-terminal graph on $n$ vertices and $m$ edges with vertex set $V(Y_n)=\{r=y_1,s=y_2,t=y_3,y_4,\cdots,y_n\}$ and edge set $E(Y_n)=\{y_iy_j|1\leq i<j\leq n\}-\{ry_4\}$.

$(3)$ Using $Z_n$ denotes the three-terminal graph on $n$ vertices and $m$ edges with vertex set $V(Z_n)=\{r=z_1,s=z_2,t=z_3,z_4,\cdots,z_n\}$ and edge set $E(Z_n)=\{z_iz_j|1\leq i<j\leq n\}-\{z_4z_5\}$.

\begin{thm}\label{th4-1}
Let $n\geq5$ and $m=\binom{n}{2}-1$ be positive integers. Then $Z_n$ is the unique uniformly most reliable graph in $\mathcal{G}_{n,m}$.
\end{thm}

\begin{pf}
By the definition of three-terminal reliability polynomial, if we can prove that there are more $rst$-subgraphs with $i$ edges in $Z_n$ than in $X_n$ and $Y_n$ for each $2\leq i\leq\binom{n}{2}-1$, then $Z_n$ is the unique uniformly most reliable graph in $\mathcal{G}_{n,m}$.

We complete this proof by construct two injective maps $f_X$ and $f_Y$, from the $rst$-subgraphs with $i$ edges in $X_n$ and $Y_n$ to the $rst$-subgraphs with $i$ edges in $Z_n$, respectively.

{\bf Construct the map $f_X$:}

Let $S$ be an $rst$-subgraph with $i$ edges in $X_n$, where $2\leq i\leq\binom{n}{2}-1$. According to whether $x_4x_5\in S$, there are two cases we need to consider.

Case 1. If $S$ does not contain the edge $x_4x_5$, then $f_X(S)=\{z_iz_j|x_ix_j\in S\}$.
The image is an $rst$-subgraph of $Z_n$ with the same number of edges as $S$. And this image does not contain the edge $rs$.

Case 2. Assume that $S$ contains the edge $x_4x_5$.

Case 2.1. If $S-\{x_4x_5\}$ is still an $rst$-subgraph, then $f_X(S)=\{z_iz_j|x_ix_j\in S\}\cup\{rs\}-\{z_4z_5\}$.
The image is an $rst$-subgraph of $Z_n$ with the same number of edges as $S$. Since this image contains the edge $rs$, it is distinct from Case 1. And $f_X(S)-\{rs\}$ is still an $rst$-subgraph.

Case 2.2. If $S-\{x_4x_5\}$ is not an $rst$-subgraph, but an $rt$-subgraph or an $st$-subgraph, then $f_X(S)=\{z_iz_j|x_ix_j\in S\}\cup\{rs\}-\{z_4z_5\}$.
The image is an $rst$-subgraph of $Z_n$ with the same number of edges as $S$. Since the image contains $rs$ and $f_X(S)-\{rs\}$ is not an $rst$-subgraph, it is distinct from the above cases. It is clear to see that in case $2.2$, $f_X(S)$ contains either the edge $st$ and an edge $rz_i$ for some $4\leq i\leq n$, or the edge $rt$ and an edge $sz_j$ for some $4\leq j\leq n$, or an edge $rz_i$ and an edge $sz_j$ for some $4\leq i,j\leq n$.

Case 2.3. Assume that $S-\{x_4x_5\}$ is not an $rst$-subgraph, or $rt$-subgraph, or $st$-subgraph.
In this case, all $rst$-subgraphs, $rt$-subgraphs and $st$-subgraphs in $S$ contain the edge $x_4x_5$. Thus, the image of the map defined by the above cases is not an $rst$-subgraph of $Z_n$. Let $S'$ be a minimal $rst$-subgraph in $S$. In $S'-\{x_4x_5\}$, there are two components containing $x_4$ and $x_5$ respectively. Targets $r$ and $s$ are in the same component, otherwise $S-\{x_4x_5\}$ is either an $rt$-subgraph or an $st$-subgraph. Without loss of generality, let $r,s$ be in the component containing $x_4$ in $S'-\{x_4x_5\}$, let $t$ be in the same component as $x_5$ in $S'-\{x_4x_5\}$.

Case 2.3.1. If $S'$ is consisted of an edge $sx_j$, a minimal $x_4x_j$-subgraph, the edge $rx_4$, the edge $x_4x_5$ and a minimal $x_5t$-subgraph, then $f_X(S)=\{z_iz_j|x_ix_j\in S\}\cup\{z_5z_j|sx_j\in S\}\cup\{rs\}-\{z_4z_5\}$.
$S$ does not have both edges $x_5x_j$ and $sx_j$, otherwise, an edge $x_5x_j$, an edge $sx_j$ and a $x_5t$-subgraph will get a $st$-subgraph that does not contain the edge $x_4x_5$ which contradicts the condition of case $2.3$. Therefore, $f_X(S)$ has the same number of edges as $S$. In $f_X(S)$, we have an $rst$-subgraph of $Z_n$ which is consisted of an edge $z_5z_j$, a $z_4z_j$-subgraph, the edge $rz_4$, a $z_5t$-subgraph and the edge $rs$. Since $f_X(S)$ contains the edge $rs$ but does not contain any edge $sz_j$ $(3\leq j\leq n)$ and $f_X(S)-\{rs\}$ is not an $rst$-subgraph, $f_X(S)$ is distinct from the above cases. The image $f_X(S)$ in this case contains the edge $rz_4$.

Case 2.3.2. If $S'$ is consisted of an edge $rx_j$, a minimal $sx_4x_j$-subgraph, the edge $x_4x_5$ and a minimal $x_5t$-subgraph, then $f_X(S)=\{z_iz_j|x_ix_j\in S\}\cup\{z_5z_j|rx_j\in S\}\cup\{rs\}-\{z_4z_5\}$.
Similarly as case $2.3.1$, $S$ does not have both edges $x_5x_j$ and $rx_j$. Therefore, $f_X(S)$ has the same number of edges as $S$. In $f_X(S)$, we have an $rst$-subgraph of $Z_n$ which is consisted of an edge $z_5z_j$, a $sz_4 z_j$-subgraph, a $z_5t$-subgraph and the edge $rs$. Since $f_X(S)$ contains edges $rs$ and $sz_j$ for some $4\leq j\leq n$ but does not contain any edge $rz_i$ $(3\leq i\leq n)$ and $f_X(S)-\{rs\}$ is not an $rst$-subgraph, $f_X(S)$ is distinct from the above cases.

Case 2.3.3. If $S'$ is consisted of the edge $rx_4$, the edge $sx_4$, the edge $x_4x_5$ and a minimal $x_5t$-subgraph, then $f_X(S)=\{z_iz_j|x_ix_j\in S\}\cup\{rz_5|sx_4\in S\}\cup\{z_5z_j|sx_j\in S,j\neq4\}\cup\{rs\}-\{z_4z_5\}$.
Similarly as case $2.3.1$, $S$ does not have both edges $x_5x_j$ and $sx_j$. Therefore, $f_X(S)$ has the same number of edges as $S$. In $f_X(S)$, we have an $rst$-subgraph of $Z_n$ is consisted of the edge $rz_4$, the edge $rz_5$, a $z_5t$-subgraph and the edge $rs$. Since $f_X(S)$ contains edges $rs$, $rz_4$ and $rz_5$ but does not contain any edge $sz_j$ $(3\leq i\leq n)$ and $f_X(S)-\{rs\}$ is not an $rst$-subgraph, $f_X(S)$ is distinct from the above cases.

Since the map $f_X$ defined on each of these cases yield $rst$-subgraph of $Z_n$ as disjoint images, the map is injective.

Because there are at least as many $rst$-subgraphs with $i$ edges in $Z_n$ as in $X_n$ for $2\leq i\leq\binom{n}{2}-1$, $Z_n$ is more reliable than $X_n$ for all $p$ $(0\leq p\leq1)$.

{\bf Construct the map $f_Y$:}

Let $S$ be an $rst$-subgraph with $i$ edges in $Y_n$, where $2\leq i\leq\binom{n}{2}-1$. According to whether $y_4y_5\in S$, there are two cases we need to consider.

Case 1. If $S$ does not contain the edge $y_4y_5$, then $f_Y(S)=\{z_iz_j|y_iy_j\in S\}$.
The image is an $rst$-subgraph of $Z_n$ with the same number of edges as $S$. And this image does not contain the edge $rz_4$.

Case 2. Assume that $S$ contains the edge $y_4y_5$.

Case 2.1. If $S-\{y_4y_5\}$ is still an $rst$-subgraph, then $f_Y(S)=\{z_iz_j|y_iy_j\in S\}\cup\{rz_4\}-\{z_4z_5\}$.
The image is an $rst$-subgraph of $Z_n$ with the same number of edges as $S$. Since the image contains the edge $rz_4$, it is distinct from Case 1. And $f_Y(S)-\{rz_4\}$ is still an $rst$-subgraph.

Case 2.2. If $S-\{y_4y_5\}$ is not an $rst$-subgraph, but is an $rsy_5$-subgraph and a $y_4t$-subgraph, or an $ry_5$-subgraph and an $sty_4$-subgraph, or an $rty_5$-subgraph and an $sy_4$-subgraph, then $f_Y(S)=\{z_iz_j|y_iy_j\in S\}\cup\{rz_4\}-\{z_4z_5\}$.
The image is an $rst$-subgraph of $Z_n$ with the same number of edges as $S$. Since the image contains $rz_4$ and $f_Y(S)-\{rz_4\}$ is not an $rst$-subgraph, it is distinct from the above cases. Since $S$ contains an edge $ry_j$ for some $2\leq j\leq n$ and $j\neq4$, the image also contains an edge $rz_j$ for some $2\leq j\leq n$ and $j\neq4$.

Case 2.3. Assume that $S-\{y_4y_5\}$ does not contain an $rst$-subgraph, or an $rsy_5$-subgraph and a $y_4t$-subgraph, or an $ry_5$-subgraph and an $sty_4$-subgraph, or an $rty_5$-subgraph and an $sy_4$-subgraph.
In this case, all $rst$-subgraphs contain the edge $y_4y_5$. And $S-\{y_4y_5\}$ is either an $rsy_4$-subgraph and a $y_5t$-subgraph, or an $rty_4$-subgraph and an $sy_5$-subgraph, or an $ry_4$-subgraph and an $sty_5$-subgraph. Therefore, the image of the map defined for the above cases is not an $rst$-subgraph of $Z_n$. Let $S'$ be a minimal $rst$-subgraph in $S$.

Case 2.3.1. If $S'$ is consisted of an edge $ry_j$, a minimal $sy_4y_j$-subgraph, the edge $y_4y_5$ and a minimal $y_5t$-subgraph, then $f_Y(S)=\{z_iz_j|y_iy_j\in S\}\cup\{z_5z_j|ry_j\in S\}\cup\{rz_4\}-\{z_4z_5\}$.
If $S$ contains both edges $y_5y_j$ and $ry_j$, then the edge $y_5y_j$, the edge $ry_j$, an $sy_4y_j$-subgraph and a $y_5t$-subgraph will yield an $rst$-subgraph without the edge $y_4y_5$, which contradicts the condition of case $2.3$. Consequently, $S$ does not have both edges $y_5y_j$ and $ry_j$ and $f_Y(S)$ has the same number of edges as $S$. In $f_Y(S)$, we have an $rst$-subgraph of $Z_n$ which is consisted of an edge $z_5z_j$, an $sz_4z_j$-subgraph, a $z_5t$-subgraph and the edge $rz_4$. Since $f_Y(S)$ contains the edge $rz_4$ but does not contain any edge $rz_j$ $(2\leq j\leq n, j\neq4)$ and $f_Y(S)-\{rz_4\}$ is not an $rst$-subgraph, $f_Y(S)$ is distinct from the above cases. In this case, each $rs$-subgraph in image $f_Y(S)$ does not contain $z_5$.

Case 2.3.2. If $S'$ is consisted of an edge $ry_j$, a minimal $ty_4y_j$-subgraph, the edge $y_4y_5$ and a minimal $sy_5$-subgraph, then $f_Y(S)=\{z_iz_j|y_iy_j\in S\}\cup\{z_5z_j|ry_j\in S\}\cup\{rz_4\}-\{z_4z_5\}$.
Similarly as case $2.3.1$, $S$ does not have both edges $y_5y_j$ and $ty_j$. Therefore, $f_Y(S)$ has the same number of edges as $S$. In $f_Y(S)$, we have an $rst$-subgraph of $Z_n$ which is consisted of an edge $z_5z_j$, a $tz_4z_j$-subgraph, a $sz_5$-subgraph and the edge $rz_4$. Since $f_Y(S)$ contains the edge $rz_4$ but does not contain any edge $rz_j$ $(2\leq j\leq n, j\neq4)$ and $f_Y(S)-\{rz_4\}$ is not an $rst$-subgraph and all $rs$-subgraph in $f_Y(S)$ contain $z_5$, $f_Y(S)$ is distinct from the above cases. In $f_Y(S)$, each $rt$-subgraph does not contain $z_5$.

Case 2.3.3. If $S'$ is consisted of an edge $ry_j$, a minimal $y_4y_j$-subgraph, the edge $y_4y_5$ and a minimal $sty_5$-subgraph, then $f_Y(S)=\{z_iz_j|y_iy_j\in S\}\cup\{z_5z_j|ry_j\in S\}\cup\{rz_4\}-\{z_4z_5\}$.
Similarly as case $2.3.1$, $S$ does not have both edges $y_5y_j$ and $ry_j$. Therefore, $f_Y(S)$ has the same number of edges as $S$.  In $f_Y(S)$, we have an $rst$-subgraph of $Z_n$ which is consisted of an edge $z_5z_j$, a $z_4z_j$-subgraph, a $stz_5$-subgraph and the edge $rz_4$. Since $f_Y(S)$ contains the edge $rz_4$ but does not contain any edge $rz_j$ $(2\leq j\leq n, j\neq4)$ and $f_Y(S)-\{rz_4\}$ is not an $rst$-subgraph and all $rs$-subgraphs and $rt$-subgraphs in $f_Y(S)$ contain $z_5$, $f_Y(S)$ is distinct from the above cases.

Since the map $f_Y$ defined on each of these cases yields $rst$-subgraph of $Z_n$ as disjoint images, the map is injective.

Because there are at least as many $rst$-subgraphs with $i$ edges in $Z_n$ as in $Y_n$ for $2\leq i\leq\binom{n}{2}-1$, $Z_n$ is more reliable than $Y_n$ for all $p$ $(0\leq p\leq1)$.

From the above argument, we conclude that $Z_n$ is the unique most reliable graph in $\mathcal{G}_{n,m}$ for all $p$ $(0\leq p\leq1)$.\qed
\end{pf}

\section{Conclusion}
The reliability of three-terminal graphs with large number of edges are investigated in this article, which fills in the blank of this research.

When the number of vertices $n$ is $4,5$ or $6$ and the number of edges $m$ is $\binom{n}{2}-2$, the uniformly most reliable graph is determined with comparisons in Example $1$ and Appendix A. When $n\geq7$ and $\binom{n}{2}-\lfloor\frac{n-3}{2}\rfloor\leq m\leq\binom{n}{2}-2$, there is no uniformly most reliable graph. However, the locally most reliable graph in $\mathcal{G}_{n,m}$ for $p$ close to $0$ or $1$ is justified by Theorems $3.1$ and $3.2$, respectively. It is worth considering whether there is a uniformly most reliable graph in the class of three-terminal graphs by deleting more edges from $K_n$ with three target vertices.

The uniformly most reliable graph in $\mathcal{G}_{n,\binom{n}{2}-1}$ is determined, which is a graph by removing an edge between non-target vertices from $K_n$ with three target vertices. This conclusion is significant in comparison with the conclusion given by Betrand \emph{et al}. {\upshape\cite{BertrandGoffGraves}, which states that the uniformly most reliable graph with $\binom{n}{2}-1$ edges for two-terminal graphs is a graph by deleting an edge between non-target vertices from $K_n$ with two target vertices. For $m=\binom{n}{2}-1$, we conjecture that the uniformly most reliable graph for $k$-terminal graphs is a graph by removing an edge between non-target vertices from $K_n$ with $k$ target vertices.

Based on these results, we found that when $m\geq\binom{n}{2}-\lfloor\frac{n-3}{2}\rfloor$, there is no uniformly most reliable three-terminal graphs in most three-terminal graph classes, then it is interesting to study the uniformly most reliable three-terminal graphs for sparse graphs or other general graphs.
In addition, it is significant to extend the results of this paper and {\upshape\cite{BertrandGoffGraves} to $k$-terminal graphs.

\begin{flushleft}
{\bf Data Availability}

No data were used to support this study.

{\bf Conflicts of Interest}

The authors declare that they have no conflicts of interest.

{\bf Acknowledgments}

This work is supported by the National Science Foundation of China (Grant Nos. 11661069, 11801296, 61603206), and the Science Found of Qinghai Province (Grant Nos. 2018-ZJ-718, 2019-ZJ-7012, 2019-ZJ-7093).

{\bf Note}

This paper has been submitted to "Mathematical Problems in Engineering".
\end{flushleft}

\newpage
\begin{appendices}
\section{Reliable polynomials for three-terminal graphs with $4<n<7$ vertices and $m=\binom{n}{2}-2$ edges}

{\bf All three-terminal graphs with $5$ vertices and $8$ edges}\\

\begin{figure}[!hbpt]
\begin{center}
\includegraphics[scale=0.5]{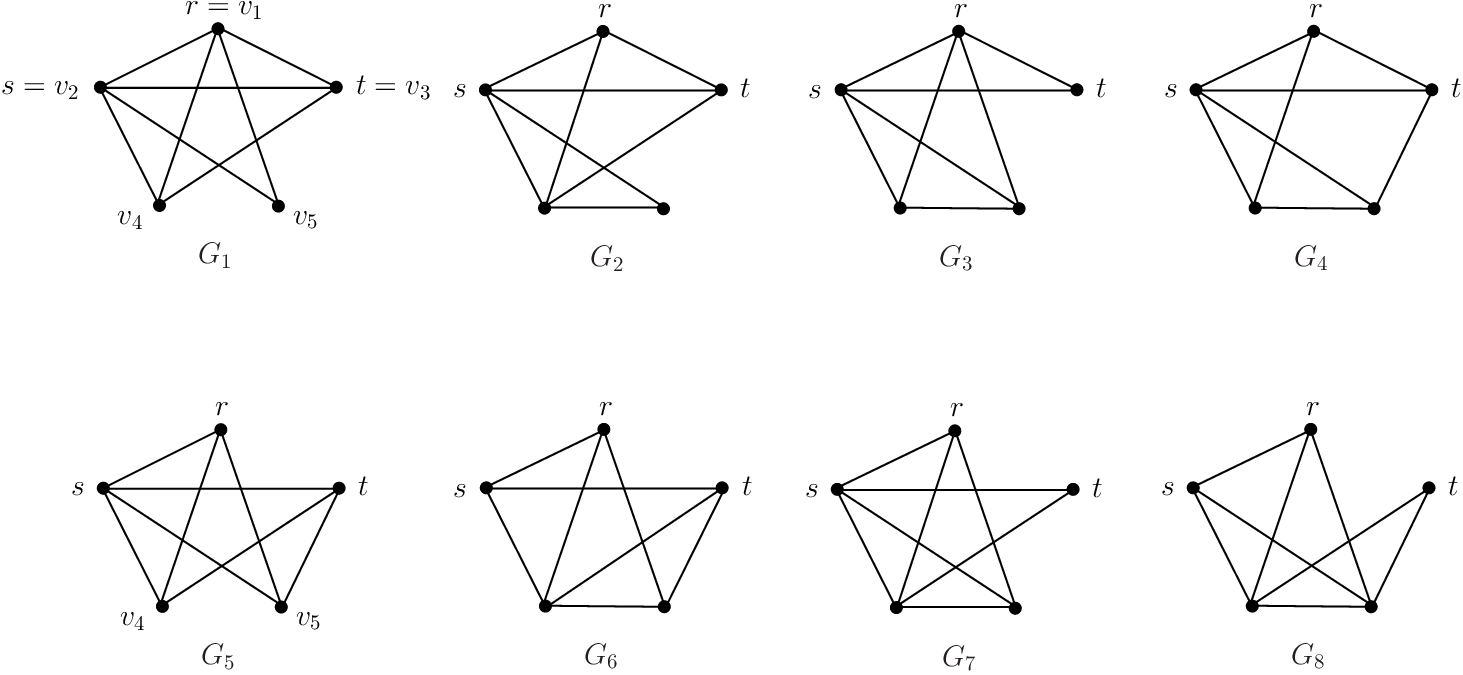} \\ [0.05cm]
Figure 6: All simple three-terminal graphs in $\mathcal{G}_{5,8}$ with three target vertices $r,s,t$.
\end{center}\label{fig20}
\end{figure}

By calculation, the reliable polynomials of these graphs are:

\begin{table}[!h]
\centering
\setlength{\abovecaptionskip}{0pt}%
\setlength{\belowcaptionskip}{6pt}%
\caption{Reliability polynomials of graphs for $R_3(G;p)$}
\setlength{\tabcolsep}{6mm}{
\begin{tabular}{l|c|c|c|c|c|c|c}
\hline  
~$R_3(G;p)$ & $N_2$ & $N_3$ & $N_4$ & $N_5$ & $N_6$ & $N_7$ & $N_8$ \\
\hline  
$R_3(G_1;p)$ & 3 & 25 & 60 & 55 & 28 & 8 & 1 \\
\hline  
$R_3(G_2;p)$ & 3 & 23 & 57 & 54 & 28 & 8 & 1 \\
\hline
$R_3(G_3;p)$ & 3 & 20 & 51 & 50 & 27 & 8 & 1 \\
\hline
$R_3(G_4;p)$ & 3 & 20 & 56 & 54 & 28 & 8 & 1\\
\hline
$R_3(G_5;p)$ & 1 & 16 & 55 & 54 & 28 & 8 & 1 \\
\hline
$R_3(G_6;p)$ & 1 & 13 & 51 & 53 & 28 & 8 & 1 \\
\hline
$R_3(G_7;p)$ & 1 & 12 & 46 & 49 & 27 & 8 & 1 \\
\hline
$R_3(G_8;p)$ & 0 & 6 & 42 & 48 & 27 & 8 & 1 \\
\hline
\end{tabular}}
\end{table}
It is clear to see that, $G_1$ is the uniformly most reliable graph in $\mathcal{G}_{5,8}$.\\

\noindent
{\bf All three-terminal graphs with $6$ vertices and $13$ edges}\\

\begin{figure}[!hbpt]
\begin{center}
\includegraphics[scale=0.8]{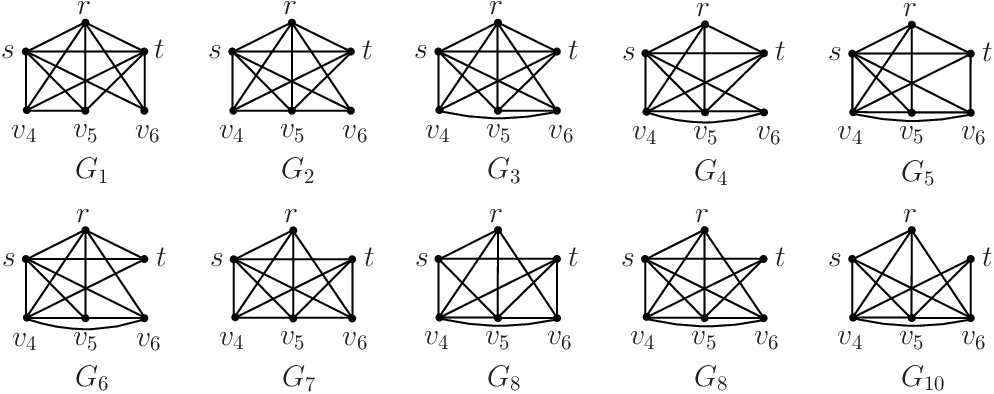} \\ [0.05cm]
Figure 7: All simple three-terminal graphs in $\mathcal{G}_{6,13}$ with three target vertices $r,s,t$.
\end{center}\label{fig20}
\end{figure}

By calculation, the reliable polynomials of these graphs are:

\begin{table}[!h]
\centering
\setlength{\abovecaptionskip}{0pt}%
\setlength{\belowcaptionskip}{6pt}%
\caption{Reliability polynomials of graphs for $R_3(G;p)$}
\setlength{\tabcolsep}{2.2mm}{
\begin{tabular}{l|c|c|c|c|c|c|c|c|c|c|c|c}
\hline  
$R_3(G;p)$ & $N_2$ & $N_3$ & $N_4$ & $N_5$ & $N_6$ & $N_7$ & $N_8$ & $N_9$ & $N_{10}$ & $N_{11}$ & $N_{12}$ & $N_{13}$\\
\hline  
$R_3(G_1;p)$ & 3 & 52 & 337 & 1017 & 1605 & 1689 & 1284 & 715 & 286 & 78 & 13 & 1 \\
\hline  
$R_3(G_2;p)$ & 3 & 47 & 304 & 955 & 1550 & 1661 & 1276 & 714 & 286 & 78 & 13 & 1 \\
\hline
$R_3(G_3;p)$ & 3 & 47 & 297 & 953 & 1552 & 1662 & 1276 & 714 & 286 & 78 & 13 & 1 \\
\hline
$R_3(G_4;p)$ & 3 & 45 & 283 & 907 & 1501 & 1634 & 1268 & 713 & 286 & 78 & 13 & 1 \\
\hline
$R_3(G_5;p)$ & 3 & 42 & 264 & 889 & 1494 & 1633 & 1268 & 713 & 286 & 78 & 13 & 1 \\
\hline
$R_3(G_6;p)$ & 3 & 42 & 259 & 849 & 1428 & 1577 & 1240 & 705 & 285 & 78 & 13 & 1 \\
\hline
$R_3(G_7;p)$ & 1 & 26 & 227 & 863 & 1486 & 1632 & 1268 & 713 & 286 & 78 & 13 & 1 \\
\hline
$R_3(G_8;p)$ & 1 & 23 & 199 & 805 & 1432 & 1604 & 1260 & 712 & 286 & 78 & 13 & 1 \\
\hline
$R_3(G_9;p)$ & 1 & 22 & 188 & 761 & 1365 & 1548 & 1232 & 704 & 285 & 78 & 13 & 1 \\
\hline
$R_3(G_{10};p)$ & 0 & 9 & 132 & 687 & 1308 & 1520 & 1224 & 703 & 285 & 78 & 13 & 1 \\
\hline
\end{tabular}}
\end{table}

It is clear to see that, $G_1$ is the uniformly most reliable graph in $\mathcal{G}_{6,13}$.
\end{appendices}

\end{document}